\documentclass[12TP,draft]{article}

\begin{document}

\begin{center}
\LARGE\noindent\textbf{ Sufficient  conditions for hamiltonian cycles in bipartite digraphs }\\

\end{center}
\begin{center}
\noindent\textbf{Samvel Kh. Darbinyan }\\

Institute for Informatics and Automation Problems, Armenian National Academy of Sciences

E-mail: samdarbin@ipia.sci.am\\
\end{center}

\textbf{Abstract}

We prove two sharp sufficient conditions for hamiltonian cycles in  balanced bipartite directed  graph. Let $D$ be a strongly connected balanced bipartite directed graph of order $2a$. Let $x,y$ be distinct vertices in $D$. $\{x,y\}$ dominates a vertex $z$ if $x\rightarrow z$ and $y\rightarrow z$; in this case, we call the pair $\{x,y\}$ dominating.  

 (i) {\it If  $a\geq 4$ and $max \{d(x), d(y)\}\geq 2a-1$ for every dominating pair of vertices $\{x,y\}$, then either $D$ is hamiltonian  or $D$ is isomorphic to one exceptional digraph of order eight.}

(ii) {\it If $a\geq 5$ and $d(x)+d(y)\geq 4a-3$ for every dominating pair of vertices $\{x,y\}$, then  $D$ is hamiltonian.} 

The first result improves a theorem of R. Wang (arXiv:1506.07949 [math.CO]), the second result, in particular, establishes a conjecture due to Bang-Jensen, Gutin and Li  (J. Graph Theory , 22(2), 1996) for strongly connected balanced bipartite digraphs of order at least ten.\\
 
\textbf{Keywords:} Digraphs, cycles, hamiltonian cycles, bipartite balanced digraph, perfect matching, longest non-hamiltonian cycles. \\

\section {Introduction} 

We consider directed graphs (digraphs) in the sense of \cite{[5]}. For convenience of the reader terminology and notations will be given in details in section 2. A digraph $D$ is hamiltonian  if it contains a cycle passing through all the vertices of $D$. For general digraphs there are many sufficient conditions for existence of hamiltonian cycles in digraphs (see,   e.g., \cite{[5]}, \cite{[7]}, \cite{[10]}, \cite{[15]}, \cite{[16]}, \cite{[21]}).
In this note, we will be concerned with the degree conditions. 

Let us recall the following well-known degree conditions (Theorems 1.1-1.4) that guarantee that a digraph is hamiltonian.
\\
\textbf{Theorem 1.1} (Nash-Williams \cite{[19]}). {\it Let $D$ be a digraph of order $n\geq 3$ such that for every vertex $x$, $d^+(x)\geq n/2$ and $d^-(x)\geq n/2$, then $D$ is Hamiltonian.}\\

 \textbf{Theorem 1.2} (Ghouila-Houri \cite{[11]}). {\it Let $D$ be a strongly connected digraph of order $n\geq 3$. If $d(x)\geq n$ for all vertices $x\in V(D)$, then $D$ is Hamiltonian.}\\

 \textbf{Theorem 1.3} (Woodall \cite{[23]}). {\it Let $D$ be a digraph of order $n\geq 3$. If $d^+(x)+d^-(y)\geq n$ for all pairs of vertices $x$ and $y$ such that there is no arc from $x$ to $y$, then $D$ is Hamiltonian.}\\
 
\textbf{Theorem 1.4} (Meyniel \cite{[18]}). {\it Let $D$ be a strongly connected digraph of order $n\geq 2$. If $d(x)+d(y)\geq 2n-1$ for all pairs of non-adjacent vertices in $D$, then $D$ is Hamiltonian.}\\

It is easy to see that Meyniel's theorem is a common generalization of Nash-Williams', Ghouila-Houri's and Woodall's theorems. For a short proof of Theorem 1.4 can be found in \cite{[9]}. \\

For bipartite digraphs, an analogue of Nash-Williams' theorem was given by Amar and Manoussakis in \cite{[3]}.
 An  analogue of Woodall's theorem was given by Manoussakis and Millis in \cite{[17]}, and  strengthened by J. Adamus and L. Adamus \cite{[1]}. The results analogous  to the above-mentioned theorems of  Ghouila-Houri and Meyniel for bipartite digraphs was given by J. Adamus, L. Adamus and A.Yeo \cite{[2]}. \\

\noindent\textbf{Theorem 1.5} (J. Adamus, L. Adamus and A.Yeo \cite{[2]}). 
 {\it Let $D$ be a balanced bipartite digraph of order $2a$, where $a\geq 2$. Then $D$ is hamiltonian provided one of the following holds:

(a) $d(u)+d(v)\geq 3a+1$ for every pair  of nonadjacent  distinct vertices $u$ and $v$ of $D$;

(b) $D$ is strongly connected and  $d(u)+d(v)\geq 3a$ for every pair 
 of nonadjacent  distinct vertices $u$ and $v$ of $D$;

(c) the minimal degree of $D$ is at least $(3a+1)/2$;

(d) $D$ is strongly connected and the minimal degree of $D$ is at least $3a/2$.}\\
 
Some sufficient conditions for the existence of a hamiltonian cycles in a bipartite tournament are described in the survey paper \cite{[13]} by Gutin.

For semicomplete bipartite digraphs a characterization was obtained independently  by Gutin \cite{[12]}  and  H\"{a}ggvist and Manoussakis \cite{[14]}.
  
Notice that each of  Theorems 1.1-1.4  imposes a degree condition on all pairs of nonadjacent vertices (or on all vertices). In the following  theorems  a degree condition requires only  for some pairs of nonadjacent vertices. \\

Let $x,y$ be distinct vertices in a digraph $D$. We say that the pair of vertices $\{x,y\}$ dominates a vertex $z$ if $x\rightarrow z$ and $y\rightarrow z$; in this case, we call the pair $\{x,y\}$ dominating. \\

\textbf{Theorem 1.6} (Bang-Jensen, Gutin, H.Li \cite{[6]}). {\it Let $D$ be a strongly connected digraph of order $n\geq 2$. Suppose that $min\{d(x),d(y)\}\geq n-1$ and  $d(x)+d(y)\geq 2n-1$ for every pair of nonadjacent vertices $x,y$ with a common in-neighbour, then $D$ is Hamiltonian.}\\

\textbf{Theorem 1.7} (Bang-Jensen, Guo, Yeo \cite{[4]}). {\it Let $D$ be a strongly connected digraph of order $n\geq 2$. Suppose that $min\{d^+(x)+d^-(y), d^-(x)+d^+(y)\}\geq n-1$ and  $d(x)+d(y)\geq 2n-1$ for every pair of nonadjacent vertices $x,y$ with a common in-neighbour or a common out-neighbour. Then $D$ is Hamiltonian.}\\

An analogue of Theorem 1.6 for bipartite balanced digraphs was given by R. Wang \cite{[22]}.

\textbf{Theorem 1.8} (R. Wang \cite{[22]}). {\it Let $D$ be a strongly connected balanced bipartite digraph of order $2a$, where $a\geq 1$. Suppose that, for every dominating pair of vertices $\{x,y\}$, either $d(x)\geq 2a-1$ and $d(y)\geq a+1$ or $d(y)\geq 2a-1$ and $d(x)\geq a+1$. Then $D$ is hamiltonian.}\\

In \cite{[6]}, Bang-Jensen, Gutin and H. Li the following conjecture was proposed.\\

\textbf{Conjecture.} {\it Let $D$ be a strongly connected digraph of order $n\geq 2$. Suppose that $d(x)+d(y)\geq 2n-1$ 
for every pair of nonadjacent vertices $x$, $y$ with a common in-neighbor. Then $D$ is Hamiltonian.}\\

The above mentioned result of Wang and the conjecture due to Bang-Jensen, Gutin and Li were  the main motivation for the present work.

Using some ideas and arguments of \cite{[22]}, in this note we prove the following  Theorems 1.9 and 1.10 (below). For $a\geq 4$ Theorem 1.9  improves the theorem of  Wang.  Theorem 1.10, in particular, establishes the conjecture of  Bang-Jensen, Gutin and H.Li in a strong form, showing it holds for strongly connected balanced bipartite digraphs of order $n\geq 10$, if $d(x)+d(y)\geq 2n-3$ for every pair of nonadjacent vertices $x$, $y$ with a common in-neighbor.\\

 \textbf{Theorem 1.9.}  {\it Let $D$ be a strongly connected balanced bipartite digraph of order $2a\geq 8$. Suppose that  $max \{d(x), d(y)\}\geq 2a-1$ for every pair of vertices $x$, $y$ with a common out-neighbor.
 If $D$ is not isomorphic to the digraph $D(8)$, then  $D$ is hamiltonian. (see, Example 1).}\\

\textbf{Theorem 1.10.}  {\it Let $D$ be a strongly connected balanced bipartite digraph of order $2a\geq 8$. Suppose that  $d(x)+ d(y)\geq 4a-3$ for every pair of vertices $x$, $y$ with a common out-neighbor. If $D$ is not isomorphic to the digraph $D(8)$, then  $D$ is hamiltonian. }\\

Of course, Theorem 1.10 is an immediate corollary of Theorem 1.9.

\section {Terminology and Notations}

Terminology  and notations not described below follow \cite{[5]}.
  In this paper we consider finite digraphs without loops and multiple arcs. 
For a digraph $D$, we denote
  by $V(D)$ the vertex set of $D$ and by  $A(D)$ the set of arcs in $D$. The order of $D$ is the number
  of its vertices. 
 The arc of a digraph $D$ directed from $x$ to $y$ is denoted by $xy$ or $x\rightarrow y$ (we also say that $x$ dominates $y$ or $y$ is out-neighbour of $x$ and $x$ is in-neighbour of $y$), and $x\leftrightarrow y$ denotes that $x\rightarrow y$ and $y\rightarrow x$ ($x\leftrightarrow y$ is called a 2-cycle). 

For disjoint subsets $A$ and  $B$ of $V(D)$  we define $A(A\rightarrow B)$ 
   as the set $\{xy\in A(D) / x\in A, y\in B\}$; $A(A,B)=A(A\rightarrow B)\cup A(B\rightarrow A)$ and 
$E(A,B)=\{xy/x\in A, y\in B\}$. 
If $x\in V(D)$ and $A=\{x\}$ we sometimes will write $x$ instead of $\{x\}$. 
$A\rightarrow B$ means that every vertex of $A$ dominates every vertex of $B$; $A\mapsto B$ means that $A\rightarrow B$ and there is no arc from $B$ to $A$.

 The out-neighborhood of a vertex $x$ is the set $N^+(x)=\{y\in V(D) / xy\in A(D)\}$ and $N^-(x)=\{y\in V(D) / yx\in A(D)\}$ is the in-neighborhood of $x$. Similarly, if $A\subseteq V(D)$, then $N^+(x,A)=\{y\in A / xy\in A(D)\}$ and $N^-(x,A)=\{y\in A / yx\in A(D)\}$. 
The out-degree of $x$ is $d^+(x)=|N^+(x)|$ and $d^-(x)=|N^-(x)|$ is the in-degree of $x$. Similarly, $d^+(x,A)=|N^+(x,A)|$ and $d^-(x,A)=|N^-(x,A)|$. The degree of the vertex $x$ in $D$ is defined as $d(x)=d^+(x)+d^-(x)$ (similarly, $d(x,A)=d^+(x,A)+d^-(x,A)$). The subdigraph of $D$ induced by a subset $A$ of $V(D)$ is denoted by $D\langle A\rangle$ or $\langle A\rangle$ for brevity. 

The path (respectively, the cycle) consisting of the distinct vertices $x_1,x_2,\ldots ,x_m$ ( $m\geq 2 $) and the arcs $x_ix_{i+1}$, $i\in [1,m-1]$  (respectively, $x_ix_{i+1}$, $i\in [1,m-1]$, and $x_mx_1$), is denoted by  $x_1x_2\cdots x_m$ (respectively, $x_1x_2\cdots x_mx_1$). 
We say that $x_1x_2\cdots x_m$ is a path from $x_1$ to $x_m$ or is an $(x_1,x_m)$-path. 
 Given a vertex $x$ of a directed path $P$ or a directed cycle $C$, we denote by $x^+$ (respectively, by $x^-$) the successor (respectively, the predecessor) of $x$ (on $P$ or $C$), and in case of ambiguity, we precise $P$ or $C$ as a subscript (that is $x^+_P$ \ldots).

A cycle that contains all the vertices of $D$  is a Hamiltonian cycle. 
 If $P$ is a path containing a subpath from $x$ to $y$ we let $P[x,y]$ denote that subpath. Similarly, if $C$ is a cycle containing vertices $x$ and $y$, $C[x,y]$ denotes the subpath of $C$ from $x$ to $y$. 
A digraph $D$ is strongly connected (or, just, strong) if there exists a path from $x$ to $y$ and a path from $y$ to $x$ for every pair of distinct vertices $x,y$.

  For an undirected graph $G$, we denote by $G^*$ the symmetric digraph obtained from $G$ by replacing every edge $xy$ with the pair $xy$, $yx$ of arcs. 
  Two distinct vertices $x$ and $y$ are adjacent if $xy\in A(D)$ or $yx\in A(D) $ (or both). For integers $a$ and $b$, $a\leq b$, let $[a,b]$  denote the set of
all the integers which are not less than $a$ and are not greater than
$b$. 

Let $H$ be a non-trivial proper subset of vertices  of a digraph $D$. An $(x,y)$-path $P$ is a $H$-bypass if $|V(P)|\geq 3$, $x\not=y$ and $V(P)\cap H=\{x,y\}$. 

A cycle factor in $D$ is a collection of vertex-disjoint cycles $C_1, C_2, \ldots , C_l$ such that $V(C_1)\cup V(C_2)\cup \ldots \cup V(C_l)= V(D)$. 
A digraph $D$ is called a bipartite digraph if there exists a partition $X$, $Y$ of $V(D)$ into two partite sets such that every arc of $D$ has its end-vertices in different partite sets. 
It is called balanced if $|X|=|Y|$. A matching from $X$ to $Y$ is a independent set of arcs with origin in $X$ and terminus in $Y$ (A  set of arcs with no common end-vertices is called independent). If $D$ is balanced, one says that such a matching is perfect if it consists of precisely $|X|$ arcs.

Let $D$ be a balanced bipartite digraph of order $2a$, where $a\geq 2$. For integer $k$, we say that
 $D$ satisfies condition $B_k$ when $max\{d(x),d(y)\}\geq 2a-2+k$ for every dominating pair of vertices $\{x,y\}$.

The underlying graph of a digraph $D$ is denoted by $UG(D)$, it contains an edge $xy$ if $x\rightarrow y$ or $y\rightarrow x$ (or both).  

\section {Examples}

In this section we present some examples of balanced bipartite digraphs which we will use in the next sections to show that the conditions of our results (the lemmas and the theorems) are sharp.\\

\textbf{Example 1.} Let $D(8)$ be a  bipartite digraph   with partite sets $X=\{x_0,x_1,x_2,x_3\}$ and 
$Y=\{y_0,y_1,y_2,y_3\}$, and let $D(8)$ contains the arcs $y_0x_1$, $y_1x_0$, $x_2y_3$, $x_3y_2$ and all the arcs of the following 2-cycles: 
$x_i\leftrightarrow y_i$, $i\in [0,3]$, $y_0\leftrightarrow x_2$, $y_0\leftrightarrow x_3$, $y_1\leftrightarrow x_2$ and  $y_1\leftrightarrow x_3$, and contains no other arcs. 

In $D(8)$ we have 
$$
d(x_2)=d(x_3)=d(y_0)=d(y_1)=7 \quad \hbox{and} \quad d(x_0)=d(x_1)=d(y_2)=d(y_3)=d(x_3)=3,
$$
and  the dominating pairs of vertices are: $\{y_0,y_1\}$, $\{y_0,y_2\}$,$\{y_0,y_3\}$,$\{y_1,y_2\}$, $\{y_1,y_3\}$, $\{x_0,x_2\}$,
$\{x_0,x_3\}$,
$\{x_1,x_2\}$, $\{x_1,x_3\}$ and $\{x_2,x_3\}$. Note that every dominating pair satisfies condition $B_1$. Since $x_0y_0x_3y_2x_2$ $y_1x_0$ is a cycle in $D(8)$, 
it is not difficult to check that $D(8)$ is strong. 

Observe that $D(8)$ is not hamiltonian. Indeed, if $C_8$ is a hamiltonian cycle in $D(8)$, then $C_8$ would be contain the arcs $x_1y_1$ and $x_0y_0$ and therefore, the path $x_1y_1x_0y_0$ or the path $x_0y_0x_1y_1$, which is impossible since $N^-(x_0)=N^-(x_1)=\{y_0,y_1\}$.\\

\textbf{Example 2.} Let $D(6)$ be a   bipartite digraph with partite sets $X=\{x_1,x_2,x_3\}$ and 
$Y=\{y_1,y_2,y_3\}$. Let $D(6)$ contains all the arcs of the following 2-cycles  $x_i\leftrightarrow y_i$, $x_1\leftrightarrow y_3$, $x_2\leftrightarrow y_3$ and the arcs  $x_2 y_1$, $x_1y_2$;  and $D(6)$ contains no other arcs. 

Notice that $d(x_1)=d(x_2)=5$, $d(y_1)=d(y_2)=3$, $d(x_3)=2$ and $d(y_3)=6$. The dominating pairs in $D(6)$ are $\{x_1,x_2\}$, $\{x_1,x_3\}$, $\{x_2,x_3\}$, $\{y_1,y_3\}$ and $\{y_2,y_3\}$ ($\{y_1,y_2\}$ is not dominating pair). It is easy to check that $D(6)$ is strong and satisfies condition $B_1$, but $UG(D(6))$ is not 2-connected.\\

\textbf{Example 3.} Let $H(6)$ be a  bipartite digraph  with partite sets $X=\{x,y,z\}$ and 
$Y=\{u,v,w\}$. Assume that $H(6)$  contains all the  arcs of the sets: $E(\{x,y\}, u)$; $E(\{u,v,w\}, x)$; $E(v, \{y,z\})$; $E(z,\{v,w\})$ and the arc $u\rightarrow z$; and contains no other arcs. 

Observe that $C_4:=xuzwx$ is a cycle in $H(6)$. The digraph $H(6)$ is strong, $d(u)=d(x)=d(v)=4$ and the dominating pairs in $H(6)$ only are
$\{x,y\}$,  $\{u,v\}$,  $\{u,w\}$ and $\{v,w\}$. Notice that $H(6)$ satisfies condition $B_0$, but contains no perfect matching from $X$ to $Y$ since for $S=\{x,y\}$ we have $N^+(S)=\{u\}$. In particular,  $H(6)$ is not hamiltonian.\\

\textbf{Example 4.} Let $D$ be a balanced bipartite digraph of order $2a\geq 8$ with partite sets $X=A\cup B\cup \{z\}$ and 
$Y=C\cup \{u,v\}$, where the subsets $A$ and $B$ are nonempty, $ A\cap B=\emptyset$, $z\notin A\cup B$ and $u, v\notin C$. 
Let  $D$  satisfies the following conditions: 

(i) the subdigraph $\langle A\cup B\cup C\cup \{z\}\rangle$ is a complete  bipartite digraph with partite sets $A\cup B\cup \{z\}$ and $C$;

(ii)   $z\rightarrow u$ and $z\leftrightarrow v$;

(iii) $N^+(u)=A$, $N^+(v)=B$; and $D$ contains no other arcs.

It is not difficult to check that $D$ is strong, $d(x)=2a-3$ for all $x\in A\cup B$, $d(y)=2a$ for all $y\in C$, 
$d(z)=2a-1$, $d(u)=|A|+1$ and $d(v)=|B|+2$. It is not difficult to check that $max\{d(b),d(c)\}\geq 2a-3$ for every dominating pair of vertices $\{b,c\}$ (i.e, $D$ satisfies condition $B_{-1}$).

 Since for $S=A\cup B$ we have $N^+(S)=C$ and $a-1=|S|>|N^+(S)|=a-2$,  by K\"{o}ning-Hall theorem $D$ contains no perfect matching from $X$ to $Y$.\\

\textbf{Example 5.} Let $D$ be  bipartite digraph of order $2a\geq 8$ with partite sets
 $X=\{x_0,x_1,x_2,$ $\ldots , x_{a-1}\}$ and $Y=\{y_0,y_1,y_2,\ldots , y_{a-1}\}$, where $\langle(X\cup Y)\setminus \{x_0,y_0\}\rangle$ is a bipartite complete digraph with partite sets $X\setminus \{x_0\}$ and $Y\setminus \{y_0\}$. $D$ contains also the following arcs $x_0y_0, y_0x_0$, $y_1x_0$ and $x_0y_1$, and $D$ contains no other arcs.

It is not difficult to check that $D$ is strong and  satisfies condition $B_0$, but $UG(D)$ is not 2-connected.\\ 

\textbf{Example 6.} Let $H'(6)$ be a bipartite digraph  with partite sets
 $X=\{x_0, x_1, x_2\}$ and $Y=\{y_0, y_1, y_2\}$. $H'(6)$  the arcs $y_0 x_1$, $y_1 x_2$, $x_0y_2$ and all the arcs of the following  2-cycles $x_i\leftrightarrow y_i$, $i\in [0,2]$,  $x_0\leftrightarrow y_1$ and $x_1\leftrightarrow y_2$:  $H'(6)$ contains no other arcs.

It is not difficult to check that $H'(6)$ is strong and  satisfies condition $B_1$, but $H'(6)$ is not hamiltonian.

\section { Preliminaries }

 \textbf{Bypass Lemma} (Lemma 3.17, Bondy \cite{[8]}). {\it Let $D$ be a strongly connected  nonseparable (i.e., $UG(D)$ is 2-connected) digraph, and let $H$ be a nontrivial proper subdigraph of $D$. Then $D$ contains a $H$-bypass.}\\

{\it Remark.} One can prove Bypass Lemma using the proof of Theorem 5.4.2 \cite{[5]}.\\

Now we prove a series of lemmas.\\

\textbf{Lemma 4.1.} {\it Let $D$ be a strong balanced bipartite digraph of order $2a\geq 4$, with partite sets $X$ and $Y$. If $ d(x)+ d(y)\geq 2a+3$ for every dominating pair of vertices $\{x,y\}$, then 

(i) $UG(D)$ is 2-connected;

(ii) if $C$ is a cycle of length $m$, $2\leq m\leq 2a-2$, then $D$ contains a $C$-bypass. }

\textbf{Proof of Lemma 4.1}. (i). Suppose, on the contrary, that $UG(D)$ is not 2-connected. Then $V(D)=A\cup B\cup \{u\}$, where $A$ and $B$ are nonempty subsets of $V(D)$, $A\cap B=\emptyset$, $u\notin A\cup B$ and there are no arcs between $A$ and $B$. Since $D$ is strong, it follows that there are two vertices $x\in A$ and $y\in B$ such that $\{x,y\}\rightarrow u$, i.e., $\{x,y\}$ is a dominating pair. Without loss of generality, assume that $x,y\in X$. Then $u\in Y$. Put $Y_1=A\cap Y$ and $Y_2=B\cap Y$. It is easy to see that $d(x)\leq 2+2|Y_1|$ and $d(y)\leq 2+2|Y_2|$. Therefore,
$$ 
d(x)+d(y)\leq 4+2(|Y_1|+|Y_2|)\leq 2a+2,
$$
which is a contradiction since $\{x,y\}$ is a dominating pair. 

(ii) The claim of Lemma 4.1(ii)  immediately follows from Lemma 4.1(i) and Bypass Lemma.
Lemma 4.1 is proved. \fbox \\\\

Note that Lemma 4.1 is not needed for the proof of  Theorem 1.9.\\

\textbf{Lemma 4.2.} {\it Let $D$ be a strong balanced bipartite digraph of order $2a\geq 8$ with partite sets $X$ and $Y$. If $D$ satisfies condition $B_1$, then

(i) $UG(D)$ is 2-connected;

(ii) if $C$ is a cycle of length $m$, $2\leq m\leq 2a-2$, then $D$ contains a $C$-bypass.}

\textbf{Proof of Lemma 4.2}. (i). Suppose, on the contrary, that $D$ is strong and satisfies condition $B_1$ but $UG(D)$ is not 2-connected. Then $V(D)= E\cup F\cup \{u\}$, where $E\cap F=\emptyset$, $u\notin E\cup F$ and there are no arcs between $E$ and $F$. Since $D$ is strong, it follows that there are two vertices $x\in E$ and $y\in F$ such that $\{x,y\}\rightarrow u$, i.e., $\{x,y\}$ is a dominating pair. By condition $B_1$, $max\{d(x),d(y)\}\geq 2a-1$.  
Without loss of generality, we assume that $x,y\in X$ and $d(x)\geq 2a-1$. Then $u\in Y$. 
From $d(x)\geq 2a-1$ and $A(E,F)=\emptyset$ it follows that $E$ contains $a-1$ vertices of $Y$, i.e., $Y\cap F=\emptyset$. 
Then, since $a\geq 4$, there exist two distinct vertices of $Y\cap E$ say $y_1,y_2$, such that $\{y_1,y_2\}\rightarrow x$, i.e., $\{y_1,y_2\}$ is a dominating pair.
Since $d(y, \{y_1,y_2\})=0$, we have $max\{d(y_1),d(y_2)\}\leq 2a-2$, which contradicts condition $B_1$. This proves that $UG(D)$ is 2-connected.

(ii). The second claims of the lemma immediately follows from the first claim and Bypass Lemma. Lemma 4.2 is proved. \fbox \\\\ 

The digraph $D(6)$ (Example 2) shows that the bound on order of $D$ in Lemma 4.2 is sharp. 

The digraph of Example 5 shows that for any $a\geq 4$ in Lemma 4.2
we cannot replace condition $B_0$ instead of $B_1$.\\

\textbf{Lemma 4.3.} {\it Let $D$ be a strong balanced bipartite digraph of order $2a\geq 8$ with partite sets $X$ and $Y$. If $D$ satisfies condition $B_0$, then $D$ contains a perfect matching from $X$ to $Y$ and a perfect matching from $Y$ to $X$. Moreover, $D$ contains a cycle factor.}

\textbf{Proof of Lemma 4.3.} By K\"{o}ning-Hall theorem (see, e.g., \cite{[5]}) to show that $D$ contains a perfect matching from $X$ to $Y$, it suffices to show that $|N^+(S)|\geq |S|$ for every set $S\subseteq X$. Let $S\subseteq X$. If $|S|=1$ or $|S|=a$, then $|N^+(S)|\geq |S|$ since $D$ is strong. Assume that $2\leq |S|\leq a-1$. We claim that $|N^+(S)|\geq |S|$. Suppose, that this is not the case, i.e., $|N^+(S)|\leq |S|-1\leq a-2$. From this and strongly connectedness of  $D$ it follows that there are two vertices $x,y\in S$ and a vertex $z\in N^+(S)$ such that $\{x,y\}\rightarrow z$, i.e., $\{x,y\}$ is a dominating pair. Hence, by condition $B_0$,
$max\{d(x),d(y)\}\geq 2a-2$. 
Without loss of generality, we assume that $d(x)\geq 2a-2$. It is easy to see that
$$
2a-2\leq d(x)\leq 2|N^+(S)|+a -|N^+(S)|=a+|N^+(S)|.
$$
Therefore, $|N^+(S)|\geq a-2$. Thus, $|N^+(S)|=a-2$ and $|S|=a-1$ since $|N^+(S)|\leq a-2$. Now it is easy to see that $d(x)=2a-2$, and hence, $\{u,v\}\rightarrow x$, where $\{u,v\}=Y\setminus N^+(S)$. By  condition $B_0$, $max\{d(u),d(v)\}\geq 2a-2$. Without loss of generality, we assume that $d(u)\geq 2a-2$. On the other hand,

$$
2a-2\leq d(u)\leq |S|+2(a -|S|)=2a-|S|,
$$
since $|S|=a-1\leq 2$. Thus, $a\leq 3$, a contradiction to $a\geq 4$.

Therefore, for any  $S\subseteq X$ we have, $|N^+(S)|\geq |S|$. By K\"{o}ning-Hall theorem there exists a perfect matching from $X$ to $Y$. The proof for a perfect matching in opposite direction is analogous. It is well-known that a digraph $D$ contains a cycle factor if and only if there exists both a perfect matching from $X$ to $Y$ and a perfect matching from $Y$ to $X$ ( Ore in \cite{[20]} Section 8.6 has shown that a balanced bipartite digraph $D$ with partite sets $X$ and $Y$ has a cycle factor if and only if $|N^+(S)|\geq |S|$ and $|N^+(T)|\geq |T|$ for all $S\subseteq X$ and $T\subseteq Y$). Therefore, $D$ contains a cycle factor. Lemma 4.3 is proved. \fbox \\\\

The digraph $H(6)$ (Example 3) shows that the bound on order of $D$ is sharp in Lemma 4.3.

 The digraph of Example 4 shows that in Lemma 4.3 instead of condition $B_0$ we cannot replace condition $B_{-1}$.\\

\textbf{Lemma 4.4.} {\it Let $D$ be a strong balanced bipartite digraph of order $2a\geq 8$ 
 with partite sets $X$ and $Y$. If $D$ is not a directed cycle  and  satisfies condition $B_0$, i.e.,  $max\{ d(x), d(y)\}\geq 2a-2$ for every dominating pair of vertices $\{x,y\}$, then $D$ contains a non-hamiltonian cycle of length at least 4.}  

\textbf{Proof of Lemma 4.4}. If $D$ is hamiltonian and  $D$ is not a directed cycle, then it is not difficult to show that $D$ contains a non-hamiltonian cycle of length at least 4. 
Suppose that $D$ is not hamiltonian and contains no cycle of length at least 4. By Lemma 4.3, $D$ contains a cycle factor. Let $C_1, C_2, \dots, C_t$ be a minimal cycle factor of $D$ (i.e., $t$ is as small  as possible). Then the length of every $C_i$ is equal to two and $t=a$. Let $C_i=x_iy_ix_i$, where $x_i\in X$ and $y_i\in Y$. 
Since $D$ is strong, there exists a vertex such that its in-degree is at least two, which means that there exists a dominating pair of vertices, say, $u$ and $v$. 
By condition $B_0$,
 $max\{d(u),d(v)\}\geq 2a-2$. 
Without loss of generality, assume that $u,v\in X$, $u=x_1$ and
$$
d(x_1)\geq 2a-2. \eqno (1)
$$
 Since $a\geq 4$ and (1), there exists a vertex $y\in Y\setminus \{y_1\}$ such that $y\leftrightarrow x_1$, say  $x_1\leftrightarrow y_2$. 
It is easy to see that $y_1$ and $x_2$ are not adjacent, else $D$ would be contains a cycle of length 4. Then $\{y_1,y_2\}$ is a dominating pair since $\{y_1,y_2\}\rightarrow x_1$.. Therefore, by condition $B_0$,
 
$$
max\{d(y_1),d(y_2)\}\geq 2a-2. \eqno (2)
$$

{ \it Case 1. $d(y_1)\geq 2a-2$.} 

Then $y_1$ together with every vertex $x_i$ other than $x_2$ forms a 2-cycle, since $y_1$ and $x_2$ are nonadjacent. This implies that $D$ contains a 4-cycle, since $x_1$ is adjacent to every vertex of $Y$, maybe except one.

{ \it Case 2. $d(y_1)\leq 2a-3$.}

Then by (2), $d(y_2)\geq 2a-2$. 

Assume first that $y_2$ together with some vertex of $X\setminus \{x_1,x_2\}$ forms a 2-cycle, say $y_2\leftrightarrow x_3$. 
Then it is not difficult to see that $max\{d(x_3),d(x_2)\}\geq 2a-2$ since $\{x_2,x_3\}$ is a dominating pair. Since $D$ contains no cycle of length 4, it is not difficult to check that   
$$
d(x_1, \{y_3\})=d(x_2, \{y_1,y_3\})=d(x_3, \{y_1\})=0. 
$$ 

These imply that $d(x_2)\leq 2a-4$, $d(x_3)\geq 2a-2$, $x_3\leftrightarrow y_4$ and $x_1\leftrightarrow y_4$. 
Therefore, $x_1y_4x_3y_2x_1$ is a cycle of length 4, which contradicts our assumption that $D$ contains no non-hamiltonian cycle of length at least 4.

Assume second that $d(y_2,\{x_i\})\leq 1$ for all $x_i\notin \{x_1,x_2\}$.
 Then from $d(y_2)\geq 2a-2$ it follows that $a=4$ and  
 $d(y_2,\{x_i\})= 1$ for $i=3$ and $i=4$. 

Let $d^+(y_2,\{x_3,x_4\})\geq 1$. 
Without loss of generality, we may assume that $y_2\rightarrow x_3$. 
Using the supposition that $D$ contains no non-hamiltonian cycle of length at least 4, it is not difficult to show that $d^+(y_3,\{x_1,x_2\})=0$, $x_3y_1\notin A(D)$. 
Notice that $d(x_2)\leq 2a-3$ since $d(x_2,\{y_1\})=d^-(x_2,\{y_3\})=0$. Thus we have
$$
A(\{x_3,y_3\}\rightarrow \{x_1,y_1,x_2,y_2\})=\emptyset . \eqno (3)
$$
If $y_2\rightarrow x_4$, by an argument similar to that in the proof of (3), we obtain  

$$
A(\{x_4,y_4\}\rightarrow \{x_1,y_1,x_2,y_2\}=\emptyset, 
$$
which together with (3) contradicts that $D$ is strong. 
Assume therefore that  $y_2 x_4\notin A(D)$. Then  $x_4\rightarrow y_2$, since $d(y_2,\{x_4\})= 1$. 
From $d(x_2)\leq 2a-3$, $\{x_2,x_4\}\rightarrow y_2$ and  condition $B_0$ it follows that $d(x_4)\geq 2a-2$.
 On the other hand, using the supposition that $D$ contains no non-hamiltonian cycle of length at least 4, it is easy to see that $d^-(x_4,\{y_1,y_3\})=0$. This together with $d(y_2,\{x_4\})=1$ gives $d(x_4)\leq 2a-3$, which is a contradiction. 

Let now  $d^+(y_2,\{x_3,x_4\})=0 $. Then $\{x_3,x_4\}\rightarrow y_2$. Now it is easy to see  that $A(\{x_1,x_2\}\rightarrow \{y_3,y_4\})=\emptyset$ and $d^+(y_1,\{x_3,x_4\})=0$. Therefore, $A(\{x_1,y_1,x_2,y_2\})\rightarrow \{x_3,y_3,x_4,y_4\}=\emptyset$, which contradicts that $D$ is strong.  Lemma 4.4 is proved. \fbox \\\\

Observe that each of digraphs $C_6^*$ and $P_6^*$ satisfies the conditions of Lemma 4.4, but has no cycle of length 4, where $C_6$ ($P_6$) is an undirected cycle (path)  with six vertices.\\

\section {The proof of the main result}

 \textbf{Proof of Theorem 1.9}. Suppose, on the contrary, that $D$ is not hamiltonian. 
In particular, $D$ is not isomorphic to the directed cycle of length $2a$. Let $C:=x_0y_0x_1y_1\ldots x_{m-1}y_{m-1}x_0$  be a longest cycle in $D$, where $x_i\in X$ and $y_i\in Y$ for all $i\in [0,m-1]$ (all subscripts are taken modulo $m$, i.e., $x_{m+i}=x_i$ and $y_{m+i}=y_i$ for all $i\in [0,m-1]$). By Lemma 4.4, $D$ contains a cycle of length at least 4, i.e., $m\geq 2$. 
By Lemma 4.2(ii), $D$ has a $C$-bypass. Let $P:=xu_1u_2\ldots u_sy$ be a $C$-bypass ($s\geq 1$). The length of the path $C[x,y]$ is the gap of $P$ with respect to $C$.
 Suppose also that the gap of $P$ is minimum among the gaps 
of all $C$-bypass. Since $C$ is a longest cycle in $D$, the length of $C[x,y]$ is greater than or equal to $s+1$.\\

First we prove that $s=1$. Suppose, on the contrary, that is $s\geq 2$.
 Since $C$ is a longest cycle in $D$ and $P$ has the minimum gap among the gaps of all $C$-bypass, $y^-_C$ is not adjacent to any vertex on $P[u_1,u_s]$ and $u_s$ is not adjacent to any vertex on $C[x^+_C, y^-_C]$. Hence,
$$ 
d(u_s)\leq 2a-2 \quad \hbox{and} \quad  d(y_C^-)\leq 2a-2
$$
since each of $P[u_1,u_s]$  and $C[x^+_C, y^-_C]$ contains at least one vertex from each partite set.
On the other hand, since $\{u_s, y^-_C\}$ is a dominating pair, by  condition $B_1$, we have
$$
max \{ d(u_s), d(y_C^-)\}\geq 2a-1,
$$
a contradiction. So, $s=1$.\\

Since $s=1$ and $D$ is a bipartite digraph, it follows  that $x$ and $y$ belong to the same partite set and the length of $C[x,y]$ must be even. Now assume, without loss of generality, that $x=x_0$, $y=x_r$ and $u_1:=v$.
 Denote $C':=V(C[y_0,y_{r-1}])$ and $R:=V(D)\setminus V(C)$. Next we consider two cases.\\

\textbf{Case 1.} $r\geq 2$.

Let $x$ be an arbitrary vertex of $X\cap R$. Since $C$ is a longest cycle in $D$, it is easy to see that 
$$
d^+(x,\{v\})+d^+(y_{r-1},\{x\})\leq 1 \quad \hbox{and} \quad d^+(v,\{x\})+d^+(x,\{y_0\})\leq 1. \eqno (4)
$$ 
Note that $\{v,y_{r-1}\}$ is a dominating pair. Recall that $v$ is not adjacent to any vertex of $C'$ since $C$-bypass $P$ has the minimum gap among the gaps of all $C$-bypass. So, for $v$ and for every $x_i\in X\cap C'$ we have,
$$ 
d(v)\leq 2a-2 \quad \hbox{and} \quad  d(x_i)\leq 2a-2. \eqno (5)
$$
Combining this with   condition $B_1$ we obtain
$$
d(y_{r-1})\geq 2a-1. \eqno (6)
$$
From this we conclude that

(i) the vertex $y_{r-1}$ and every vertex of $X$ are adjacent. In particular, $y_{r-1}$ and  $x$ are adjacent, i.e., $x\rightarrow y_{r-1}$ or $y_{r-1}\rightarrow x$.\\

If $x\rightarrow y_{r-1}$, then $d^-(x,\{v,y_0\})=0$ since $C$-bypass $P$ has the minimum gap among the gaps of all $C$-bypass. Hence, $d(x)\leq 2a-2$. This together with $d(x_{r-1})\leq 2a-2$ (by (5)) gives a contradiction since $\{x,x_{r-1}\}$ is a dominating pair. Assume therefore that $xy_{r-1}\notin A(D)$. Then $y_{r-1}\mapsto x$. By the arbitrariness of $x$, we may assume that $y_{r-1}\mapsto X\cap R$. Combining this with $d(y_{r-1})\geq 2a-1$ (by (6)), we obtain that $|R|=2$, i.e., the cycle $C$ has length $2a-2$ and 

(ii) the vertex $y_{r-1}$ and every vertex of $X\cap V(C)$ form a 2-cycle. In particular, $\{x_0,x_1\}\rightarrow y_{r-1}$ and $d(x_0)\geq 2a-1$, since by (5),  $d(x_1)\leq 2a-2$.\\

Therefore, any two distinct vertices of $X\cap V(C)$ form a dominating pair, which means that every vertex of $X\cap V(C)$, except for at most one vertex, has degree at least  $2a-1$. This together with the second inequality of (5) and
condition $B_1$  implies that $r=2$ and for every $x_i\in \{x_0,x_1,\ldots , x_{m-1}\}\setminus \{x_1\}$,
$$
d(x_i)\geq 2a-1,    \eqno (7)
$$ 
which also means that

(iii) the vertex $v$ and every vertex  $x_i\notin \{x_0,x_1,\ldots , x_{m-1}\}\setminus \{x_1\}$ are adjacent.\\

Using the first inequality of (4) and $x y_1\notin A(D)$ (by our assumption) we obtain that $d^+(x,\{v,y_1\})=0$. Therefore, since $D$ is strong, it follows that there is a vertex $y_l$ other than $y_1$, such that $x\rightarrow y_l$.
 Notice that $l\not= 2$, since $P$ has the minimum gap among the gaps of all $C$-bypass. If $l\leq m-1$ and $x_l\rightarrow y_0$, then $vx_2\ldots x_ly_0x_1y_1xy_l\ldots x_0v$ is a hamiltonian cycle, a contradiction. Thus we can assume that
$$
\hbox{if} \quad 3\leq l\leq m-1, \quad \hbox{then} \quad x_ly_0\notin A(D). \eqno (8)
$$
We now will consider the following subcases below.

\textbf{Subcase 1.1.} $v\rightarrow x_0$.

From the minimality of the gap $|C[x_0,x_r]|-1$ of $P$ and (iii) it follows that $v\mapsto \{x_2,x_3,\ldots , x_{m-1}\}$. This together with (7) implies that

(iv) every vertex $x_i$, other than $x_0$ and $x_1$, and every vertex $y_j$ form a 2-cycle. In particular, for all $i\in [2,m-1]$, $x_i\leftrightarrow y_0$ and $x_2\rightarrow \{y_0,y_1\}$.\\

Now using (8) and (iv), it is not difficult to show that $l=0$, i.e.,
$$
x\rightarrow y_0 \quad \hbox{and} \quad d^+(x,\{y_1,y_2,\ldots , y_{m-1}\})=0.  \eqno (9)
$$
From $y_1\rightarrow x\rightarrow y_0$ and (4) we obtain that $v$ and $x$ are not adjacent. Therefore, $d(x)\leq 2a-2$.  
This together with $d(x_1)\leq 2a-2$ (by (5)) and  condition  $B_1$ implies that $ x_1 y_0\notin A(D)$. Since $\{v,y_0\}\rightarrow x_2$ and $d(v)\leq 2a-2$ (by (5)), from  condition $B_1$ it follows that $d(y_0)\geq 2a-1$.
 This and $ x_1 y_0\notin A(D)$ give $y_0\rightarrow x_0$.
 Combining this with (iv) we obtain that $y_0\rightarrow \{x_2,x_3,\ldots x_{m-1},x_0\}$. 
Therfore, if for some $j\in [2,m-1]$, $y_j\rightarrow x_1$, then, by (9), we have that $x\rightarrow y_0$ and  $vx_2\ldots y_jx_1y_1xy_0x_{j+1}\ldots x_0v$ is a hamiltonian cycle in $D$, a contradiction. Assume therefore that 
$$
d^-(x_1,\{y_2,y_3,\ldots , y_{m-1}\})=0.
$$
The last equality together with (9) implies that $d(y_j)\leq 2a-2$ for every $j\in [2,m-1]$. Now recall that $\{v,y_2\}\rightarrow x_2$, by (iv) (i.e., $\{v,y_2\}$ is a dominating pair). But $d(y_2)\leq 2a-2$ and $d(v)\leq 2a-2$, which contradicts  condition $B_1$.

\textbf{Subcase 1.2.} $x_2\rightarrow v$ and $v x_0\notin A(D)$.

Then from the minimality of the gap $|C[x_0,x_2]|-1$ and (iii) it follows that 
$$
\{x_3,x_4,\ldots , x_{m-1},x_0\}\mapsto v.  \eqno (10)
$$
This together with (7) implies that 

(v) every vertex $x_i$, other than $x_1$ and $ x_2$, and every vertex $y_j$ form a 2-cycle for all $j\in [0,m-1]$. In particular, $y_j\rightarrow x_0$ and if $i\notin \{1,2\}$, then $x_i\leftrightarrow y_0$. 

From (8) and (v) it follows that $l=0$, i.e.,
$$
x\rightarrow y_0 \quad \hbox{and} \quad d^+(x,\{y_1,y_2,\ldots , y_{m-1})=0.  \eqno (11)
$$
Since $x\rightarrow y_0$ and  $y_1\rightarrow x$, by (4) we have that $v$ and $x$  are not adjacent. 
Therefore, $d(x)\leq 2a-2$ and $x_1 y_0\notin A(D)$ because of $d(x_1)\leq 2a-2$ (by (5)) and $x\rightarrow y_0$. 
If for some $j\in [2,m-1]$,  $y_j\rightarrow x_1$, then by (v) and  $x\rightarrow y_0$ we have that $vx_2\ldots y_jx_1y_1xy_0x_{j+1}\ldots x_0v$ is a hamiltonian cycle, which contradicts our supposition that $D$ is not hamiltonian. 
Assume therefore that $d^-(x_1,\{y_2,y_3,\ldots , y_{m-1}\})=0$. 
This together with (11) implies that
$$
d(y_j)\leq 2a-2 \quad \hbox{for all} \quad j\in [2,m-1].   \eqno (12)
$$
By (v) we have $y_j\rightarrow x_0$ for all $y_j$. 
Combining this with (12) we obtain that $m=3$, i.e., the cycle $C$ has length 6 and $a=4$.
 From $d(y_2)\leq 2a-2$ (by (12)),
$\{y_0,y_2\}\rightarrow x_0$ (by (v)) and  condition $B_1$ it follows that  $d(y_0)\geq 2a-1$. 
Then $y_0\rightarrow x$ and $y_0\leftrightarrow x_2$ since $x_1y_0\notin A(D)$.
 Using (ii), i.e., the fact that the vertex $y_1$ forms a 2-cycle with each vertex of $\{x_0,x_1,x_2\}$ it is easy to show that $x_1$ and $y_2$ are not adjacent ( for otherwise, if $x_1\rightarrow y_2$, then $x_1y_2x_0vx_2y_1xy_0x_1$ is a hamiltonian cycle,   if $y_2\rightarrow x_1$, then $y_2x_1y_1xy_0x_0vx_2y_2$ is a hamiltonian cycle). 
Thus  we have $d(v)=3$ and $d(y_2)\leq 4$ since $d(y_2,\{x,x_1\})=0$, which implies that $y_2x_2\notin A(D)$, because of $v\rightarrow x_2$, $d(v)\leq 3$ and condition $B_1$. 
Thus we have $a=4$, $D$ contains exactly the following 2-cycles and arcs: $v\leftrightarrow x_2$, 
$x_2\leftrightarrow y_1$,
$y_1\leftrightarrow x_1$, $y_1\leftrightarrow x_0$, $y_2\leftrightarrow x_0$, $y_0\leftrightarrow x_0$, $y_0\leftrightarrow x$, $y_0\leftrightarrow x_2$, $x_2 y_2$, $y_1 x$, $y_0 x_1$ and  $x_0 v$.

Now it is not difficult to see that $D$ is isomorphic to $D(8)$. (To see this, let now $X:=\{x_0,x_1,x_2, x_3\}$ and  $Y:=\{y_0,y_1,y_2, y_3\}$, where $x_0:=x$, $x_1:=x_1$, $x_2:=x_0$, $x_3:=x_2$, $y_0:=y_0$, $y_1:=y_1$, $y_2:=y_2$ and $y_3:=v$). Subcase 1.2 is considered. 

\textbf{Subcase 1.3.} $x_2v\notin A(D)$ and $vx_0\notin A(D)$.

Let the vertex $v$ and $t$ vertices of $C[x_3,x_{m-1}]$ form a 2-cycle (recall that $v\in Y$). We will consider the subcases $t\geq 1$ and $t=0$  separately.

\textbf{Subcase 1.3.1.} $t\geq 1$.

Then $m\geq 4$. Let $x_q\in C[x_3,x_{m-1}]$ be a vertex such that $v$ and $x_q$ form a 2-cycle and $q$ is minimal with these properties. 
From this, (iii) and the fact that $C$-bypass $P$ has the minimum gap among the gaps of all  $C$-bypass it follows that
$$
v\mapsto \{x_2,x_3,\ldots , x_{q-1}\} \quad \hbox{and} \quad \{x_{q+1},x_{q+2},\ldots , x_{m-1}, x_0\}\mapsto v. \eqno (13)
$$
Hence, $t=1$. From (13) and (7) it follows

(vi) every vertex $x_i\in C[x_2,x_0]\setminus \{x_q\}$ together with every vertex $y_j$ forms a 2-cycle. In particular, 
$$ 
y_j\leftrightarrow x_2 \quad \hbox{and} \quad x_i\leftrightarrow y_0      \eqno (14)
$$
for all $y_j$ and for all $x_i$ other than $x_1$ and $x_q$, respectively. So, $\{v,y_j\}\rightarrow x_2$, i.e., $\{v,y_j\}$ is a dominating pair for all $y_j$. 
From condition $B_1$ and $d(v)\leq 2a-2$ (by (5))  for all $j\in [0,m-1]$ we have,

$$ d(y_j)\geq 2a-1. \eqno (15)
$$
Notice that from $x_i\leftrightarrow y_0$ (by (14)) and (8) imply that $l=q$ or $l=0$. Recall that $x\rightarrow y_l$.

Let $l=q$, i.e., $x\rightarrow y_q$.
 By (8), $x_qy_0\notin A(D)$. This together with $d(x_q)\geq 2a-1$ (by (7)) implies that $y_0\rightarrow x_q\rightarrow y_{q-1}$ and $y_1\rightarrow x_q$. 
Recall that $x_{q-1}\leftrightarrow y_0$ by (14). Since $C$-bypass $P$ has the minimum gap among the gaps of all $C$-bypass, it follows that $y_{q-1}x\notin A(D)$ (for otherwise, the $C$-bypass $y_{q-1}\rightarrow x\rightarrow y_q$ has a gap equal to 2 which is a contradiction). 
This together with $d(y_{q-1})\geq 2a-1$ (by (15)) implies that $y_{q-1}\rightarrow x_1$, Thus,
 $vx_2\ldots y_{q-1}x_1y_1xy_q\ldots x_0y_0x_qv$ is a Hamiltonian cycle in $D$, which contradicts our initial supposition.

Let now $l=0$, i.e., 
$$
x\rightarrow y_0 \quad \hbox{and} \quad d^+(x, C[y_1,y_{m-1}])=0.    \eqno (16)
$$   
Then, by (4), the vertices $v$ and $x$ are not adjacent and hence, $d(x)\leq 2a-2$. From $d(x_1)\leq 2a-2$ (by (5)), $d(x)\leq 2a-2$ and  condition $B_1$ we have $x_1y_0\notin A(D)$. 
Since $d(y_{m-1})\geq 2a-1$ (by (15)) and $xy_{m-1}\notin A(D)$ (by (16)), it follows that $x_1\leftrightarrow y_{m-1}$. Therefore, if $y_0\rightarrow x_0$, then $vx_2\ldots y_{m-1}x_1y_1xy_0x_0v$ is a hamiltonian cycle, a contradiction.
 Assume then that $y_0x_0\notin A(D)$. Then from $vx_0\notin A(D)$ it follows that $d(x_0)\leq 2a-2$ which contradicts that $d(x_0)\geq 2a-1$ by (7).

\textbf{Subcase 1.3.2.} $t=0$, i.e., there is no $x_i$, $i\in [0,m-1]$, such that $x_i\leftrightarrow v$.

From (7) it follows that

(vii) every vertex $x_i$ other than $x_1$ and every vertex $y_j$ form a 2-cycle. In particular, for every $i\not= 1$ and every $j\in [1,m-1]$ we have
$$
 x_i\leftrightarrow y_0 \quad \hbox{and}\quad  y_j\leftrightarrow x_2.     \eqno (17)
$$
Now using (8) and  $x_i\leftrightarrow y_0$, $i\not= 1$, we obtain  
$$
 d^+(x,C[y_1,y_{m-1}])=0  \quad \hbox{and}\quad  x\rightarrow y_0.     \eqno (18)
$$
By (4), it is easy to see that $v$ and $x$ are not adjacent because of $y_1\rightarrow x$ and $x\rightarrow y_0$. Therefore, $d(x)\leq 2a-2$.
If $y_j\rightarrow x_1$ for some $y_j\notin \{y_0,y_1\}$, then, by (17), $y_0\rightarrow x_{j+1}$, and hence $vx_2\ldots y_jx_1y_1xy_0x_{j+1}\ldots x_0v$ is a hamiltonian cycle, a contradiction.
Assume therefore that 
$d^-(x_1,\{y_2,y_3,\ldots , y_{m-1}\})=0$.
This together with the first equality of (18) implies that 
$$
 d(y_j)\leq 2a-2  \quad \hbox{for all}\quad  y_j\notin \{y_0,y_1\},     \eqno (19)
$$
in particular, $d(y_2)\leq 2a-2$. From (17) we have, $y_2\rightarrow x_2$. 
 Hence, $\{v,y_2\}\rightarrow x_2$, i.e., $\{v,y_2\}$ is a dominating pair, but by (5) and (19) we have 
$max\{d(v),d(y_2)\}\leq 2a-2$, 
which contradicts  condition $B_1$, and completes the discussion  of Case 1.\\

\textbf{Case 2.} $r=1$.

Note that $\{v,y_0\}$ is a dominating pair. By condition $B_1$,
$max\{d(v),d(y_0)\}\geq 2a-1$.
Assume, without loss of generality, that $d(v)\geq 2a-1$, which implies that

(vii) the vertex $v$ and every vertex of $X$ are adjacent.\\

Note that $X\cap R\not=\emptyset$ and consider the following three subcases.

\textbf{Subcase 2.1.} There exists a vertex $u\in X\cap R$ such that $v\leftrightarrow u$.

Using our supposition that $C$ is a longest cycle in $D$, it is not difficult to show that the following two claims are true.

 \textbf{Claim 1.} For any $x_i\in V(C)$ if $x_i\rightarrow v$, then $u y_i\notin A(D)$;
 if $v\rightarrow x_i$, then $y_{i-1} u\notin A(D)$.\\

\textbf{Claim 2.} If there exists $x_i\in V(C)$ such that $x_i\rightarrow v\rightarrow x_{i+1}$, then $u$ and $y_i$ are not adjacent.\\

Now we will prove the following claim.

\textbf{Claim 3.} If there exists $x_i\in V(C)$ such that $x_i\leftrightarrow v$, then (a) $x_{i+1}\mapsto v$ or (b) $v\mapsto x_{i-1}$ is impossible.

\textbf{Proof of Claim 3.} (a). Suppose, on the contrary, that  $x_i\leftrightarrow v$ and  $x_{i+1}\mapsto v$.
Combining this with $d(v)\geq 2a-1$ (by our assumption) we obtain,
$$
 v\leftrightarrow x \quad \hbox{for every}\quad  x\in X\setminus \{x_{i+1}\}.     \eqno (20)
$$
This with Claim 2 implies 

(viii) if $v\leftrightarrow z$ where $z\in X\cap R$, then $z$ and every vertex of $(Y\cap V(C))\setminus \{y_i\}$ are not adjacent. In particular,  the following holds
$$
d(y_j)\leq 2a-2 \quad \hbox{and} \quad d(z)\leq 2a-2    \eqno (21)
$$
for any $y_j$ other than $y_i$.\\

If $|X\cap R|\geq 2$, then by (20) and (21) there are two distinct vertices in $X\cap R$, say $x,z$, such that
$x\leftrightarrow v$, $z\leftrightarrow v$ and 
$max\{d(x),d(z)\}\leq 2a-2$, which contradicts  condition $B_1$.
Assume therefore that $|X\cap R|=1$.
 Then the cycle $C$ has length $2a-2$, i.e., $m=a-1\geq 3$. Since $u\leftrightarrow v$, $x_{i+1}\mapsto v$ and $d(u)\leq 2a-2$ (by (21)), from  condition $B_1$ it follows that $d(x_{i+1})\geq 2a-1$. 
This together with $x_{i+1}\mapsto v$ implies that $x_{i+1}$ and every vertex of $Y\setminus \{v\}$ form a 2-cycle, i.e., any two distinct vertices of $Y\cap V(C)$ form a dominating pair. On the other hand, since $m\geq 3$ and (21), for two distinct vertices of $(Y\cap V(C))\setminus \{y_i\}$, say $y_s$ and $y_k$, we have $max\{d(y_s),d(y_k)\}\leq 2a-2$, which contradicts  condition $B_1$. Claim 3(a) is proved.

(b). Suppose, on the contrary, that  $x_i\leftrightarrow v$ and  $v\mapsto x_{i-1}$.
Combining this with $d(v)\geq 2a-1$ (by our assumption) we obtain that $v$ and every vertex of $X\setminus \{x_{i-1}\}$ form a 2-cycle, i.e., 
$$
 v\leftrightarrow x \quad \hbox{for any}\quad  x\in X\setminus \{x_{i-1}\}.     \eqno (22)
$$
This with Claim 2 implies

(ix) every vertex of $X\cap R$  and every vertex of $(Y\cap V(C))\setminus \{y_{i-1}\}$ are not adjacent. Hence, for all $x\in X\cap R$  and  for all $y_j$ other than $y_{i-1}$ the following holds
$$
d(x)\leq 2a-2 \quad \hbox{and} \quad d(y_j)\leq 2a-2.     \eqno (23)
$$
 If $|X\cap R|\geq 2$, then by (22) and (23) for any two distinct vertices $x,z$ of $X\cap R$, we have $\{x,z\}\rightarrow v$ and  $max\{d(x),d(z)\}\leq 2a-2$, which contradicts  condition $B_1$.
Assume therefore that $|X\cap R|=1$.  
Then the cycle $C$ has length $2a-2$, i.e., $m=a-1\geq 3$. 
Since $x_i\leftrightarrow v$ and $v\leftrightarrow u$, from  condition $B_1$ and the first inequality of (23) it follows that $d(x_i)\geq 2a-1$.
 Therefore, $x_i$ and every vertex of $Y\cap V(C)$, maybe except one, form a 2-cycle. 
This together with second inequality of (23) implies that $m=3$, $y\rightarrow x_i$ for some $y\in (Y\cap V(C))\setminus \{y_{i-1}\}$ and $d(y_{i-1})\geq 2a-1$. 
It follows that $u\rightarrow y_{i-1}$ since $y_{i-1}u\notin A(D)$. 
Thus we have, $\{x_{i-1},u\}\rightarrow y_{i-1}$ (i.e.,
 $\{x_{i-1},u\}$ is a dominating pair) and  $d(x_{i-1})\geq 2a-1$ because of $d(u)\leq 2a-2$ by (23).
 Now $d(x_{i-1})\geq 2a-1$ and $v\mapsto x_{i-1}$ imply that $x_{i-1}$ and every vertex of $Y\cap V(C)$ form a 2-cycle (i.e., there are two distinct vertices of $(Y\cap V(C))\setminus \{y_{i-1}\}$, say $y_s$ and $y_k$, such that $\{y_s,y_k\}\rightarrow x_{i-1}$) 
 which is a contradiction,
because of the second inequality of (23). Claim 3 is proved. \fbox \\\\
 
Now we can finish the proof of Theorem 1.9.

By Claim 3 and $d(v)\geq 2a-1$ (by our assumption), we have that $v$ and every vertex of $X\cap V(C)$ form a 2-cycle.
Combining this with Claim 2 we obtain that  $u$ and every vertex of $V(C)$ are not adjacent, which in turn implies that $d(u)\leq 2a-2$ and  $d(y_j)\leq 2a-2$ for all $y_j$.
  Using $\{x_i,u\}\rightarrow v$, $d(u)\leq 2a-2$ and  condition $B_1$ we obtain that $d(x_i)\geq 2a-1$ for all $x_i$. Since $d(y_j)\leq 2a-2$ for all $y_j$, using condition $B_1$ 
we conclude that no two distinct vertices of $\{y_0,y_1,\ldots , ,y_{m-1}\}$ form a dominating pair. 
In particular, $y_i x_i\notin A(D)$ for all $y_i$. 
From this and $d(x_i)\geq 2a-1$ it follows that every $x_i$ together with every vertex of $Y\setminus \{y_i\}$ forms a 2-cycle.

If $|R|=2$, then $m\geq 3$ since $a\geq 4$, and $\{y_0,y_2\}\rightarrow x_1$, but $max\{d(y_0),d(y_2)\}\leq 2a-2$, which is a contradiction.

Therefore assume that $|R|\geq 4$. As noted above,  $x_i\leftrightarrow y$ for all $y\in (Y\cap R)\setminus \{v\}$. Notice that $\{y_{i-1},y\}\rightarrow x_i$, i.e., $\{y_{i-1},y\}$ is a dominating pair. From $d(y_0)\leq 2a-2$ and  condition $B_1$ we have, $d(y)\geq 2a-1$. Therefore, $y\rightarrow u$ or $u\rightarrow  y$. In both cases it is not difficult to show that $D$ contains a  cycle of length $2m+2$ which contradicts that $C$ is longest cycle in $D$. 
The discussion of Subcase 2.1 is completed.\\

By Subcase 2.1, we can assume that there is no 2-cycle between $v$ and any vertex of $X\cap R$.
From this and $d(v)\geq 2a-1$ it follows that $|R|=2$. Let $X\cap R=\{u\}$.\\

\textbf{Subcase 2.2.} $u\mapsto v$.

From  $d(v)\geq 2a-1$ and $u\mapsto v$ it follows that $v$ and every vertex of $X\setminus \{u\}$ form a 2-cycle. 
Since $D$ is strong, there exists $y_i$ such that $y_i\rightarrow u$. Hence, $y_iuvx_{i+1}\ldots x_iy_i$ is a hamiltonian cycle, a contradiction.\\

\textbf{Subcase 2.3.} $v\mapsto u$.

Similar to Subcase 2.2, we may obtain a contradiction. The theorem is proved. \fbox \\\\

The digraph $H'(6)$ (Example 6) and its converse  digraph show that the bound on order of $D$ in Theorem 1.9 is sharp.

The digraph of Example 5 shows that in Theorem 1.9 instead of condition $B_1$ we cannot replace condition $B_0$.\\

From Theorem 1.9 it immediately follows Theorem 1.10 and the following corollary.\\

\textbf{Corollary 5.1.} (Wang \cite{[22]}). {\it Let $D$ be a strongly connected balanced bipartite digraph of order $2a$, where $a\geq 5$. Suppose that, for every dominating pair of vertices $\{x,y\}$, either $d(x)\geq 2a-1$ and $d(y)\geq a+1$ or $d(y)\geq 2a-1$ and $d(x)\geq a+1$. Then $D$ is hamiltonian.}\\

\section {Concluding remarks}

A balanced bipartite digraph of order $2a$ is even pancyclic if it contains a cycle of length $2k$ for every $k$, $2\leq k \leq a$.

Motivated by the Bondy's  "metaconjecture" (see, e.g., \cite{[8]} p. 88), it is natural to set the  following problem:

{\it Characterize those balanced bipartite digraphs which satisfy condition $B_1$ but are not even pancyclic}.\\

We have proved the following theorem.
 
\textbf{Theorem 6.1.} {\it Let $D$ be a strongly connected balanced bipartite digraph of order $2a\geq 8$ other than the directed cycle of length $2a$. If $ max\{d(x), d(y)\}\geq 2a-1$ for every dominating pair of vertices $\{x,y\}$,  then either $D$ contains  cycles of all even lengths less than or equal to $2a$ or $D$ is isomorphic to the digraph $D(8)$ (Example 1).}\\

\end{document}